\documentclass[11pt]{amsart}
\usepackage{amsfonts}
\usepackage{amsmath}
\usepackage{amssymb}
\usepackage{amscd}
\usepackage[mathscr]{eucal}
\usepackage{url}

\newcommand{\ph}[2]{{\left({#1}\right)}_{#2}}

\newcommand{\Bigph}[2]{{\Bigl({#1}\Bigr)}_{#2}}

\newcommand{\gf}[1]{\Gamma{\left({#1}\right)}}
\newcommand{\biggf}[1]{\Gamma{\bigl({#1}\bigr)}}
\newcommand{\Biggf}[1]{\Gamma{\Bigl({#1}\Bigr)}}

\renewcommand*{\bar}{\overline}
\newcommand{\gfp}[1]{\Gamma_{p}{\left({#1}\right)}}

\newcommand{\bin}[2]{\left({\genfrac{}{}{0pt}{}{#1}{#2}}\right)}

\newcommand{\fac}[1]{{#1}\text{!}}

\theoremstyle{plain}
\newtheorem{theorem}{Theorem}[section]
\newtheorem{lemma}[theorem]{Lemma}
\newtheorem{prop}[theorem]{Proposition}

\theoremstyle{definition}

\newtheorem{conj}[theorem]{Conjecture}

\numberwithin{equation}{section}

\begin{document}

\title[A $p$-adic analogue of a formula of Ramanujan]{A $p$-adic analogue of a formula of Ramanujan}
\author{Dermot M\lowercase{c}Carthy and Robert Osburn}  

\address{School of Mathematical Sciences, University College Dublin, Belfield, Dublin 4, Ireland}

\email{dermot.mc-carthy@ucdconnect.ie}

\email{robert.osburn@ucd.ie}


\date{August 24, 2007}

\begin{abstract}
During his lifetime, Ramanujan provided many formulae relating binomial sums to special values of the Gamma function. Based on numerical computations, Van Hamme recently conjectured $p$-adic analogues to such formulae. Using a combination of ordinary and Gaussian hypergeometric series, we prove one of these conjectures.
\end{abstract}

\maketitle

\section{Introduction}
In Ramanujan's second letter to Hardy dated February 27, 1913, the following formula appears:

\begin{equation} \label{Ramanujan}
1 - 5 {\left({ \frac{1}{2} }\right)}^5 + 9 {\left({ \frac{1 \cdot 3}{2 \cdot 4} }\right)}^5
- 13 {\left({ \frac{1\cdot 3 \cdot 5}{2 \cdot 4 \cdot 6} }\right)}^5 + \dotsm
= \frac{2}{\gf{\frac{3}{4}}^4}
\end{equation}

\noindent where $\gf{\cdot}$ is the Gamma function. This result was proved in 1924 by Hardy \cite{H} and a further proof was given by Watson \cite{Wa} in 1931. Note that (\ref{Ramanujan}) can be expressed as
\begin{equation*}
\sum^{\infty}_{k=0} (4k+1) {\bin{-\frac{1}{2}}{k}}^5
= \frac{2}{{\gf{\frac{3}{4}}}^4} \; .
\end{equation*}

\noindent Other formulae of this type include 

\begin{equation} \label{other}
\sum_{k=0}^{\infty} (-1)^k \frac{6k+1}{4^k} {\bin{-\frac{1}{2}}{k}}^3 = \frac{4}{\pi}=\frac{4}{\gf{\frac{1}{2}}^2},
\end{equation}

\noindent which is Entry $20$, page 352 of  \cite{B1}. It is interesting to note that a proof of (\ref{other}) was not found until 1987 \cite{borwein}. 

Recently, Van Hamme \cite{VH} studied a $p$-adic analogue of (\ref{Ramanujan}). Namely, he truncated the left-hand side and replaced the Gamma function with the $p$-adic Gamma function. Based on numerical computations, he posed the following.

\begin{conj}\label{VanHamme}
Let $p$ be an odd prime. Then
\begin{equation*}
\sum^\frac{p-1}{2}_{k=0}(4k+1){\bin{-\frac{1}{2}}{k}}^5 \equiv 
\left\{ \begin{array}{ll} -\frac{p}{\gfp{\frac{3}{4}}^4} \pmod {p^3} & \qquad \textup{if} \quad p\equiv1 \pmod 4\vspace{.05in} \\
\qquad 0 \quad \pmod {p^3} & \qquad \textup{if} \quad p\equiv3 \pmod 4 \end{array} \right.
\end{equation*}
where $\gfp{\cdot}$ is the $p$-adic Gamma function.
\end{conj}

\noindent The purpose of this paper is to prove the following.

\begin{theorem}\label{The_Theorem}
Conjecture \ref{VanHamme} is true.
\end{theorem}

The paper is organized as follows. In Section 2 we recall some properties of the Gamma function, ordinary hypergeometric series, the $p$-adic Gamma function and Gaussian hypergeometric series. The proof of Theorem \ref{The_Theorem} is then given in Section 3. We would like to point out that Conjecture \ref{VanHamme} is but one indication of the interplay between ordinary hypergeometric series and Gaussian hypergeometric series. Further evidence of this interplay can be found in \cite{gs} and \cite{rouse}. Finally, van Hamme states $12$ other conjectures relating truncated hypergeometric series to values of the $p$-adic Gamma function. This includes a conjectural $p$-adic analogue of (\ref{other}) which states

\begin{equation*}
\sum_{k=0}^{\frac{p-1}{2}} (-1)^k \frac{6k+1}{4^k} {\bin{-\frac{1}{2}}{k}}^3 \equiv -\frac{p}{\gfp{\frac{1}{2}}^2} \pmod {p^4}.
\end{equation*} 

\noindent These conjectures were motivated experimentally and as van Hamme states that ``we have no real explanation for our observations", it might be worthwhile to determine whether these congruences arise from considering some appropriate algebraic surfaces.

\section{Preliminaries}

We briefly discuss some preliminaries which we will need in Section 3. For further details see \cite{AAR}, \cite{Ba}, \cite{BEW}, \cite{Ko} or \cite{Mu}. Recall that for all complex numbers $x \neq 0,-1,-2, \dotsc$, the Gamma function $\gf{x}$ is defined by

\begin{equation*}
\gf{x} := \lim_{k \rightarrow \infty} \frac{{k!} \: k^{x-1}}{\ph{x}{k}} \; 
\end{equation*}

\noindent where $\ph{a}{0}:=1$ and $\ph{a}{n} := a(a+1)(a+2)\dotsm(a+n-1)$ for positive integers $n$. The Gamma function satisfies the reflection formula

\begin{equation} \label{reflect}
\gf{x}\gf{1-x} = \frac{\pi}{\sin{\pi x}}. 
\end{equation}

We also recall that the hypergeometric series $_pF_q$ is defined by

\begin{equation}\label{GHS_def}
{_pF_q} \left[ \begin{array}{ccccc} a_1, & a_2, & a_3, & \dotsc, & a_p \vspace{.05in}\\
\phantom{a_1} & b_1, & b_2, & \dotsc, & b_q \end{array}
\Big| \; z \right]
:=\sum^{\infty}_{n=0}
\frac{\ph{a_1}{n} \ph{a_2}{n} \ph{a_3}{n} \dotsm \ph{a_p}{n}}
{\ph{b_1}{n} \ph{b_2}{n} \dotsm \ph{b_q}{n}}
\; \frac{z^n}{{n!}}
\end{equation}
where $a_i$, $b_i$ and $z$ are complex numbers, with none of the $b_i$ being negative integers or zero, and, $p$ and $q$ are positive integers. Note that the series terminates if some $a_j$ is a negative integer. In \cite{W}, Whipple studied properties of {\it well-poised} series where $p=q+1$, $z=\pm 1$, and $a_1 +1 = a_2 + b_1 = a_3 +b_2 = \dotsm = a_p + b_q$. One such transformation property of the well-poised series (see (6.3), page 252 in \cite{W}) is 

\begin{multline}\label{Whipple_Formula}\raisetag{29pt}
{_6F_5} \left[ \begin{array}{cccccc} a, & 1+\frac{1}{2}a, & c, & d, & e, & f \vspace{.05in}\\
\phantom{a} & \frac{1}{2}a, & 1+a-c, & 1+a-d, & 1+a-e, & 1+a-f \end{array}
\Big| \; -1 \right]\\
=\frac{\gf{1+a-e}\gf{1+a-f}}{\gf{1+a}\gf{1+a-e-f}}
\; {_3F_2} \left[ \begin{array}{ccc} 1+a-c-d, & e, & f \vspace{.05in}\\
\phantom{1+a-c-d} & 1+a-c, & 1+a-d \end{array}
\Big| \; 1 \right] . \; \; \quad
\end{multline}
 
\noindent This is Entry 31, Chapter 10 in Ramanujan's second notebook (see page 41 of \cite{B}). Watson's proof of (\ref{Ramanujan}) is a specialization of (\ref{Whipple_Formula}) combined with Dixon's theorem \cite{dixon}.

Let $p$ be an odd prime.  For $n \in \mathbb{N}$, we define the $p$-adic Gamma function as
\begin{align*}
\gfp{n} &:= {(-1)}^n \prod_{\substack{j<n\\p \nmid j}} j \\
\intertext{and extend to all $x \in\mathbb{Z}_p$ by setting}
\gfp{x} &:= \lim_{n \rightarrow x} \gfp{n}
\end{align*}
where $n$ runs through any sequence of positive integers $p$-adically approaching $x$ and 
$\gfp{0}:=1$. 
This limit exists, is independent of how $n$ approaches $x$ and determines a continuous function
on $\mathbb{Z}_p$.

In \cite{G}, Greene introduced the notion of general hypergeometric series over finite fields or \emph{Gaussian hypergeometric series}.  These series are analogous to classical hypergeometric series and have played an important role in relation to the number of points over $\mathbb{F}_{p}$ of Calabi-Yau threefolds \cite{AO}, traces of Hecke operators \cite{fop}, formulas for Ramanujan's $\tau$-function \cite{p}, and the number of points on a family of elliptic curves \cite{jf}. 

We now introduce two definitions. Let $\mathbb{F}_{p}$ denote the finite field with $p$ elements. We extend the domain of all characters $\chi$ of $\mathbb{F}^{*}_{p}$ to $\mathbb{F}_{p}$, by defining $\chi(0):=0$ for $\chi \neq \epsilon_p$ and $\epsilon_p(0):=1$, where $\epsilon_p$ is the trivial character mod $p$. The first definition is the finite field analogue of the binomial coefficient. For characters $A$ and $B$ of $\mathbb{F}_{p}$, define $\bin{A}{B}$ by
\begin{equation*}
\binom{A}{B} := \frac{B(-1)}{p} J(A, \bar{B})
\end{equation*}
where $J(\chi, \lambda)$  denotes the Jacobi sum for $\chi$ and $\lambda$ characters of $\mathbb{F}_{p}$. The second definition is the finite field analogue of ordinary hypergeometric series. For characters $A_0,A_1,\dotsc, A_n$ and $B_1, \dotsc, B_n$ of $\mathbb{F}_{p}$ and 
$x \in \mathbb{F}_{p}$, define the \textit{Gaussian hypergeometric series} by
\begin{equation*}
{_{n+1}F_n} {\left( \begin{array}{cccc} A_0, & A_1, & \dotsc, & A_n \\
\phantom{A_0} & B_1, & \dotsc, & B_n \end{array}
\Big| \; x \right)}_{p}
:= \frac{p}{p-1} \sum_{\chi} \binom{A_0 \chi}{\chi} \binom{A_1 \chi}{B_1 \chi}
\dotsm \binom{A_n \chi}{B_n \chi} \chi(x)
\end{equation*}
where the summation is over all characters $\chi$ on $\mathbb{F}_{p}$.

In \cite{OS}, the case where $A_i = \phi_p$, the quadratic character, for all $i$ and $B_j= \epsilon_p$ for all $j$ is examined and is denoted ${_{n+1}F_n}(x)$ for brevity.  By \cite{G}, $p^{n} {}_{n+1}F_{n}(x) \in \mathbb{Z}$. Before stating the main result of \cite{OS}, we recall that for $i$, $n \in \mathbb{N}$, generalized harmonic sums, ${H}^{(i)}_{n}$, are defined by
\begin{equation*}
{H}^{(i)}_{n}:= \sum^{n}_{j=1} \frac{1}{j^i}
\end{equation*}
and ${H}^{(i)}_{0}:=0$.
For $p$ an odd prime, $\lambda \in \mathbb{F}_{p}$, $n \in \mathbb{Z}^{+}$, we now define the quantities
\begin{multline}\label{X_Def}
X(p,\lambda,n) := \phi_p(\lambda)
\sum^{\frac{p-1}{2}}_{j=0} {\bin{\frac{p-1}{2}+j}{j}}^{l}
{\bin{\frac{p-1}{2}}{j}}^{l} {{(-1)}^{jl}} {\lambda^{-j}}
\left[1+2(n+1)j \left({H_{\frac{p-1}{2}+j}^{(1)}  }\right. \right. \\
\left. \left. {- H_{j}^{(1)}}\right)
 +\frac{{\left(n+1\right)}^2}{2} j^2 {\left({H_{\frac{p-1}{2}+j}^{(1)} - H_{j}^{(1)}}\right)}^2
-\frac{\left(n+1\right)}{2} j^2 \left({H_{\frac{p-1}{2}+j}^{(2)} - H_{j}^{(2)}}\right)\right],
\end{multline}
\begin{multline}\label{Y_Def}
Y(p,\lambda,n) := \phi_p(\lambda)
\sum^{\frac{p-1}{2}}_{j=0} {\bin{\frac{p-1}{2}+j}{j}}^{l}
{\bin{\frac{p-1}{2}}{j}}^{l} {{(-1)}^{jl}} {\lambda^{-jp}}
\left[1+(n+1)j \left(H_{\frac{p-1}{2}+j}^{(1)} \right. \right. \\
\left. \left. - H_{j}^{(1)} \right)
 -\frac{n+1}{2}j \left(H_{\frac{p-1}{2}+j}^{(1)} - H_{\frac{p-1}{2}-j}^{(1)}\right)\right],
\end{multline}
and
\begin{equation}\label{Z_Def}
Z(p, \lambda, n) := \phi_p(\lambda) \sum^{\frac{p-1}{2}}_{j=0}
{\bin{2j}{j}}^{2l} 16^{-jl} {\lambda^{-jp^2}},
\end{equation}
where $l=\frac{n+1}{2}$.
The main result in \cite{OS} provides an expression for ${_{n+1}F_{n}}$ modulo $p^3$. Precisely, we have
\begin{theorem}\label{OS_Theorem}
Let $p$ be an odd prime, $\lambda \in \mathbb{F}_{p}$, and $n \geq 2$ be an integer. Then
\begin{equation*}
-p^n {_{n+1}F_{n}} (\lambda) \equiv
{\left({-\phi_p(-1)}\right)}^{n+1}
\left[p^2 X(p,\lambda,n) + pY(p,\lambda,n) + Z(p,\lambda,n)\right] \pmod {p^3} .
\end{equation*}
\end{theorem}

\section{Proof of Theorem \ref{The_Theorem}}

\noindent {\it{Proof of Theorem \ref{The_Theorem}.}} By Proposition 4.2 in \cite{mort2} and Corollary 5 in \cite{VH}, we have that
\begin{equation*}
p^{3}  { _{3}F_{2}} (1)
= \left\{ \begin{array}{ll} -\frac{p}{{\gfp{\frac{3}{4}}}^4} \pmod {p^3} & \qquad \textup{if} \quad p\equiv1 \pmod 4 \vspace{.05in} \\
\qquad 0 \quad \pmod {p^3} & \qquad \textup{if} \quad p\equiv3 \pmod 4 \; .\end{array} \right.
\end{equation*}
Thus, by Theorem \ref{OS_Theorem} it suffices to prove
\begin{equation}\label{XYZ_Thm}
\sum^\frac{p-1}{2}_{k=0}(4k+1){\bin{-\frac{1}{2}}{k}}^5 \equiv
\phi_p(-1) \left[ p^3 X(p,1,2) + p^2 Y(p,1,2) + p Z(p,1,2)\right] \pmod {p^3}
\end{equation}
where the quantities $X(p,\lambda,n), Y(p,\lambda,n)$ and $Z(p,\lambda,n)$ are defined by (\ref{X_Def}),  (\ref{Y_Def}) and (\ref{Z_Def}) respectively.
We first show, via the following lemmas, that the terms involving $Y(p,1,2)$ and $X(p,1,2)$ in (\ref{XYZ_Thm}) vanish modulo $p^3$.

%%%%%%%%%%%%%%%% Start Lemma 1 - Y(p,1,2) %%%%%%%%%%%%%%%%%%%%%
\begin{lemma}\label{lemma1}
Let $p$ be an odd prime. Then
\[Y(p,1,2) \equiv 0 \pmod {p} \; .\]
\end{lemma}

\begin{proof}
Substituting $\lambda=1$ and $n=2$ in equation (\ref{Y_Def}), we get
\begin{multline*}
Y(p,1,2) = \sum^{\frac{p-1}{2}}_{j=0} {\bin{\frac{p-1}{2}+j}{j}}^{\frac{3}{2}}
{\bin{\frac{p-1}{2}}{j}}^{\frac{3}{2}} {{(-1)}^{\frac{3}{2}j}}
\left[1+3j \left(H_{\frac{p-1}{2}+j}^{(1)} - H_{j}^{(1)}\right) \right. \\
\left. -\frac{3}{2}j \left(H_{\frac{p-1}{2}+j}^{(1)} - H_{\frac{p-1}{2}-j}^{(1)}\right)\right] \; .
\end{multline*}

\noindent
Noting that $\binom{u+k}{k} ={(-1)}^k\binom{-1-u}{k}$, we get 
\begin{align}\label{Change_to_Minus_Half}
{\bin{\frac{p-1}{2}+j}{j}} {\bin{\frac{p-1}{2}}{j}}
={(-1)}^j \; {\bin{-\frac{1}{2}-\frac{p}{2}}{j}} {\bin{-\frac{1}{2}+\frac{p}{2}}{j}}
\equiv  {(-1)}^j {\bin{-\frac{1}{2}}{j}}^2 \pmod {p^2} \; .
\end{align}

\noindent
Also,
\begin{align*}
H_{\frac{p-1}{2}+j}^{(1)} - H_{\frac{p-1}{2}-j}^{(1)}
&=\frac{1}{\frac{p-1}{2}-j+1}+\frac{1}{\frac{p-1}{2}-j+2} + \dotsm + \frac{1}{\frac{p-1}{2}}
+\frac{1}{\frac{p+1}{2}} + \dotsm + \frac{1}{\frac{p-1}{2}+j}\\
&=\sum^{j-1}_{r=0} \; \frac{1}{\frac{p-1}{2}-r} \; + \; \frac{1}{\frac{p+1}{2}+r}\\
&=\sum^{j-1}_{r=0} \; \frac{4p}{p^2-{(2r+1)}^2}\\
&\equiv 0 \pmod {p} \; .
\end{align*}

\noindent
So we need only show 
\begin{equation}\label{lemma1_1}
\sum^{\frac{p-1}{2}}_{j=0} {\bin{-\frac{1}{2}}{j}}^3 {{(-1)}^{3j}}
\left[1+3j \left(H_{\frac{p-1}{2}+j}^{(1)} - H_{j}^{(1)}\right)\right]
\equiv 0 \pmod {p} \; .
\end{equation}

\noindent
For $j\geq1$, note that
\begin{align}\label{Change_to_Pochhammer}
{\bin{-\frac{1}{2}}{j}} {{(-1)}^{j}}\equiv \frac{{(j+1)}_{\frac{p-1}{2}}}{\fac{\left(\frac{p-1}{2}\right)}} \pmod {p} \; .
\end{align}

\noindent
As $\gcd \left( {\left(\frac{p-1}{2}\right)!}^3 \:,p \right)=1$, it now suffices to show
\begin{equation}\label{lemma1_2}
{{\left(\frac{p-1}{2}\right)}!}^3+\sum^{\frac{p-1}{2}}_{j=1} {{(j+1)}_{\frac{p-1}{2}}^3}
\left[1+3j \left(H_{\frac{p-1}{2}+j}^{(1)} - H_{j}^{(1)}\right)\right]
\equiv 0 \pmod {p} \; .
\end{equation}

\noindent
We now use an argument similar to that in Section 4 of \cite{K}. Let
\begin{equation}\label{P_Def}
P(z):=\frac{d}{dz} \left[z{(z+1)}_{\frac{p-1}{2}}^3\right]
=\sum^{\frac{3p-3}{2}}_{k=0}a_kz^k
\end{equation}
for some integers $a_k$. By a computation, we have
\begin{align*}
P(z) = {{(z+1)}_{\frac{p-1}{2}}^3}
\left[{ 1 + 3z \left({ H_{\frac{p-1}{2}+z}^{(1)} - H_{z}^{(1)} }\right) }\right] \; .
\end{align*}

\noindent
Combining this with (\ref{lemma1_2}), it is enough to show that
\begin{equation}\label{lemma1_3}
{{\left(\frac{p-1}{2}\right)}{!}}^3 + \sum^{\frac{p-1}{2}}_{j=1} P(j) \equiv 0 \pmod {p} \; .
\end{equation}

\noindent
Note that, for $\frac{p-1}{2} < j < p$, ${(j+1)}_{\frac{p-1}{2}}$ is divisible by $p$ and ${H_{\frac{p-1}{2}+j}^{(i)} - H_{j}^{(i)}} \in \frac{1}{p^i} \mathbb{Z}_p$, so that $P(j) \equiv 0 \pmod {p}$ for such $j$. Hence (\ref{lemma1_3}) will hold if we can show
\begin{equation}\label{lemma1_4}
{{\left(\frac{p-1}{2}\right)}{\text{!}}}^3 + \sum^{p-1}_{j=1} P(j) \equiv 0 \pmod {p} \; .
\end{equation}

\noindent
We now recall the following elementary fact about exponential sums. For a positive integer $k$, we have
\begin{equation}\label{exp_sums}
\sum^{p-1}_{j=1} j^k\equiv
\begin{cases}
-1 \pmod {p}& \text{if $(p-1) \vert k$} \; ,\\
\phantom{-}0 \pmod {p}& \text{otherwise} \; .
\end{cases}
\end{equation}

\noindent
By (\ref{P_Def}), (\ref{exp_sums}) and the fact that $\frac{3p-3}{2} < {2p-2}$, we see that
\begin{equation*}
\begin{split}
\sum^{p-1}_{j=1} P(j)
&=\sum^{p-1}_{j=1} \sum^{\frac{3p-3}{2}}_{k=0}a_k j^k\\
&=\sum^{\frac{3p-3}{2}}_{k=0}a_k \sum^{p-1}_{j=1} j^k\\
&\equiv -a_0 - a_{p-1} \pmod {p} \; .
\end{split}
\end{equation*}

\noindent
Additionally, by (\ref{P_Def})
\begin{equation*}
{{(z+1)}_{\frac{p-1}{2}}^3}
= \dotsm + \frac{a_{p-1}}{p}z^{p-1} + \dotsm \qquad .
\end{equation*}

\noindent As ${{(z+1)}_{\frac{p-1}{2}}^3}$ has integer coefficients, $p$ divides ${a_{p-1}}$.
Hence $a_{p-1} \equiv 0 \pmod{p}$. One can also check that
\begin{equation*}
a_0 = {{\left(\frac{p-1}{2}\right)}{!}}^3 \; .
\end{equation*}

\noindent Thus
\begin{equation*}
\sum^{p-1}_{j=1} P(j) \equiv -{{\left(\frac{p-1}{2}\right)}{!}}^3 \pmod{p}
\end{equation*}
and (\ref{lemma1_4}) holds. This proves the result.
\end{proof}

%%%%%%%%%%%% End Lemma 1 - Y(p,1,2) %%%%%%%%%%%%%%%%%%%%%%%%%

Now we would like to show that ord$_p(X(p,1,2)) \geq 0$ which ensures that the term
involving $X(p,1,2)$ in equation (\ref{XYZ_Thm}) vanishes modulo $p^3$.  In fact, in the following lemma,
we show that ord$_p(X(p,1,2)) \geq 1$.

%%%%%%%%%%%%%%%% Start Lemma 2 - X(p,1,2) %%%%%%%%%%%%%%%%%%%%%
\begin{lemma}\label{lemma2}
Let $p$ be an odd prime. Then
\[X(p,1,2) \equiv 0 \pmod {p} \; .\]
\end{lemma}

\begin{proof}
Substituting $\lambda=1$ and $n=2$ in equation (\ref{X_Def})
and applying (\ref{Change_to_Minus_Half}) and (\ref{lemma1_1}) yields
\begin{multline*}
X(p,1,2) \equiv \sum^{\frac{p-1}{2}}_{j=0} {\bin{-\frac{1}{2}}{j}}^3 {{(-1)}^{3j}}
\left[3j \left({H_{\frac{p-1}{2}+j}^{(1)} - H_{j}^{(1)}}\right)\right.\\
\left. +\frac{9}{2} j^2 {\left({H_{\frac{p-1}{2}+j}^{(1)} - H_{j}^{(1)}}\right)}^2
-\frac{3}{2} j^2 \left({H_{\frac{p-1}{2}+j}^{(2)} - H_{j}^{(2)}}\right)\right]
\pmod {p} \; .
\end{multline*}

\noindent
By (\ref{Change_to_Pochhammer}) and as $\gcd \left( {{{\left({\frac{p-1}{2}}\right)}!}^3,p} \right)=1$, it suffices to prove that
\begin{multline}\label{lemma2_1}
\sum^{\frac{p-1}{2}}_{j=1} {{(j+1)}_{\frac{p-1}{2}}^3}
\left[3j \left({H_{\frac{p-1}{2}+j}^{(1)} - H_{j}^{(1)}}\right)\right.
+\frac{9}{2} j^2 {\left({H_{\frac{p-1}{2}+j}^{(1)} - H_{j}^{(1)}}\right)}^2\\
\left. -\frac{3}{2} j^2 \left({H_{\frac{p-1}{2}+j}^{(2)} - H_{j}^{(2)}}\right)\right]
\equiv 0 \pmod {p} \; .
\end{multline}

\noindent Similar to the proof of Lemma \ref{lemma1}, we now let
\begin{equation}\label{Q_Def}
Q(z):= \frac{z}{2} \frac{d^2}{{dz}^2} \left[{z{{(z+1)}_{\frac{p-1}{2}}^3}}\right]
=\sum^{\frac{3p-3}{2}}_{k=0}a_kz^k
\end{equation} 
for some integers $a_k$. One can check that it now suffices to show
\begin{equation}\label{lemma2_3}
\sum^{p-1}_{j=1} Q(j) \equiv 0 \pmod {p} \; .
\end{equation}

\noindent
By (\ref{exp_sums}), (\ref{Q_Def}) and the fact that $\frac{3p-3}{2} < {2p-2}$, we have 
\begin{equation*}
\begin{split}
\sum^{p-1}_{j=1} Q(j)
&=\sum^{p-1}_{j=1} \sum^{\frac{3p-3}{2}}_{k=0}a_k j^k\\
&=\sum^{\frac{3p-3}{2}}_{k=0}a_k \sum^{p-1}_{j=1} j^k\\
&\equiv - a_{p-1} \pmod {p} \; .
\end{split}
\end{equation*}
Here we have used that $a_0=0$ as $z \vert Q(z)$.
One can check that
\begin{equation*}
{{(z+1)}_{\frac{p-1}{2}}^3}
= \dotsm + \frac{2 a_{p-1}}{p(p-1)} z^{p-1} + \dotsm \quad .
\end{equation*}

\noindent As ${{(z+1)}_{\frac{p-1}{2}}^3}$ has integer coefficients, $p$ divides ${a_{p-1}}$.
Hence $a_{p-1} \equiv 0 \pmod{p}$. Thus (\ref{lemma2_3}) holds and the result is proven
\end{proof}
%%%%%%%%%%%% End Lemma 2 - Y(p,1,2) %%%%%%%%%%%%%%%%%%%%%%%%%

Via (\ref{XYZ_Thm}), Lemma~\ref{lemma1} and Lemma~\ref{lemma2}, the proof of Theorem \ref{The_Theorem} is complete  on proving the following Proposition.

%%%%%%%%%%%%%%% Begin Lemma 3 - Z(p,1,2) %%%%%%%%%%%%%%%%%%%%%

\begin{prop}\label{prop3}
Let $p$ be an odd prime. Then
\begin{equation*}
\sum^{\frac{p-1}{2}}_{k=0}(4k+1){\bin{-\frac{1}{2}}{k}}^5 \equiv
\phi_p(-1) \; p \; Z(p,1,2) \pmod {p^3} \; .
\end{equation*}
\end{prop}

\begin{proof}
Substituting $\lambda=1$ and $n=2$ in equation (\ref{Z_Def}), we get
\begin{equation}\label{lemma3_1}
Z(p,1,2) = \sum^{\frac{p-1}{2}}_{j=0} {\bin{2j}{j}}^3 {16}^{-\frac{3}{2}j} \: .
\end{equation}

\noindent
Noting that
\begin{align*}
\bin{2j}{j}
&=2^{2j} {(-1)}^j \bin{-\frac{1}{2}}{j},
\end{align*}
it suffices to prove
\begin{equation}\label{lemma3_2}
\sum^{\frac{p-1}{2}}_{k=0}(4k+1){\bin{-\frac{1}{2}}{k}}^5 \equiv
\phi_p(-1) \: p
\left[{ \sum^{\frac{p-1}{2}}_{j=0} {(-1)}^j {\bin{-\frac{1}{2}}{j}}^3 }\right] \pmod {p^3} \; .
\end{equation}

%\noindent
%We now consider the transformation formula due to Whipple (see (6.3), page 252 in \cite{W}), which is a well known result in the theory of generalised hypergeometric series (for further details see \cite{Ba}). Namely,
%\begin{multline*}
%{_6F_5} \left[ \begin{array}{cccccc} a, & 1+\frac{1}{2}a, & c, & d, & e, & f, \\
%\phantom{a} & \frac{1}{2}a, & 1+a-c, & 1+a-d, & 1+a-e, & 1+a-f \end{array}
%\Big| \; -1 \right]\\
%=\frac{\gf{1+a-e}\gf{1+a-f}}{\gf{1+a}\gf{1+a-e-f}}
%\; {_3F_2} \left[ \begin{array}{ccc} 1+a-c-d, & e, & f \\
%\phantom{1+a-c-d} & 1+a-c, & 1+a-d \end{array}
%\Big| \; 1 \right] \; .
%\end{multline*}

%\noindent
%Note that if $f$ is a negative integer both the series in the above equation will terminate 
%after $(-f+1)$ terms in their respective sums. 

\noindent Letting $a=\frac{1}{2}$, $c=\frac{1}{2}+i\frac{p}{2}$, $d=\frac{1}{2}-i\frac{p}{2}$, 
$e=\frac{1}{2}+\frac{p}{2}$ and $f=\frac{1}{2}-\frac{p}{2}$ in (\ref{Whipple_Formula}), we get
\begin{multline}\label{Whipple}
{_6F_5} \left[ \begin{array}{cccccc} \frac{1}{2}, & \frac{5}{4}, & \frac{1}{2}+i\frac{p}{2} &
\frac{1}{2}-i\frac{p}{2}, & \frac{1}{2}+\frac{p}{2}, & \frac{1}{2}-\frac{p}{2} \vspace{.05in} \\
\phantom{\frac{1}{2}} & \frac{1}{4}, & 1-i\frac{p}{2}, & 1+i\frac{p}{2}, &
1-\frac{p}{2}, & 1+\frac{p}{2} \end{array}
\Big| \; -1 \right]\\
=\frac{\biggf{1-\frac{p}{2}}\biggf{1+\frac{p}{2}}}
{\biggf{\frac{3}{2}}\biggf{\frac{1}{2}}}
\; {_3F_2} \left[ \begin{array}{ccc} \frac{1}{2} & \frac{1}{2} + \frac{p}{2}, & \frac{1}{2} - \frac{p}{2} \vspace{.05in}\\
\phantom{\frac{1}{2}} & 1-i\frac{p}{2}, & 1+i\frac{p}{2} \end{array}
\Big| \; 1 \right] \; .
\end{multline}

\noindent
By (\ref{GHS_def}),
\begin{multline}\label{Whipple_LHS}
{_6F_5} \left[ \begin{array}{cccccc} \frac{1}{2}, & \frac{5}{4}, & \frac{1}{2}+i\frac{p}{2} &
\frac{1}{2}-i\frac{p}{2}, & \frac{1}{2}+\frac{p}{2}, & \frac{1}{2}-\frac{p}{2} \vspace{.05in}\\
\phantom{\frac{1}{2}} & \frac{1}{4}, & 1-i\frac{p}{2}, & 1+i\frac{p}{2}, &
1-\frac{p}{2}, & 1+\frac{p}{2} \end{array}
\Big| \; -1 \right]\\*
= \sum^{\frac{p-1}{2}}_{k=0}
\frac{\Bigph{\frac{1}{2}}{k} \Bigph{\frac{5}{4}}{k}
\Bigph{\frac{1}{2}+i\frac{p}{2}}{k} \Bigph{\frac{1}{2}-i\frac{p}{2}}{k}
\Bigph{\frac{1}{2}+\frac{p}{2}}{k} \Bigph{\frac{1}{2}-\frac{p}{2}}{k}}
{\Bigph{\frac{1}{4}}{k}
\Bigph{1-i\frac{p}{2}}{k} \Bigph{1+i\frac{p}{2}}{k}
\Bigph{1-\frac{p}{2}}{k} \Bigph{1+\frac{p}{2}}{k}}
\; \frac{(-1)^k}{{k}!} \; .
\end{multline}

\noindent
Now, 
\begin{align}\label{LHS_Sub1}
\ph{\frac{1}{2}}{k} \frac{(-1)^k}{{k}!}
=\bin{-\frac{1}{2}}{k} \: ,
\end{align}
\begin{align}\label{LHS_Sub2}
\frac{\Bigph{\frac{5}{4}}{k}}{\Bigph{\frac{1}{4}}{k}}
=4k+1 \; ,
\end{align}
and
\begin{align}\label{LHS_Sub3}
\frac
{\Bigph{\frac{1}{2}+i\frac{p}{2}}{k}  \Bigph{\frac{1}{2}-i\frac{p}{2}}{k}
 \Bigph{\frac{1}{2}+\frac{p}{2}}{k} \Bigph{\frac{1}{2}-\frac{p}{2}}{k}} 
{\Bigph{1-i\frac{p}{2}}{k}  \Bigph{1+i\frac{p}{2}}{k}
\Bigph{1-\frac{p}{2}}{k} \Bigph{1+\frac{p}{2}}{k}}
\equiv {\bin{-\frac{1}{2}}{k}}^4  \pmod {p^4} \; .
\end{align}

\noindent
Therefore, substituting (\ref{LHS_Sub1}), (\ref{LHS_Sub2}) and (\ref{LHS_Sub3})
into equation (\ref{Whipple_LHS}), we get
\begin{multline}\label{LHS_Final}
{_6F_5} \left[ \begin{array}{cccccc} \frac{1}{2}, & \frac{5}{4}, & \frac{1}{2}+i\frac{p}{2}, &
\frac{1}{2}-i\frac{p}{2}, & \frac{1}{2}+\frac{p}{2}, & \frac{1}{2}-\frac{p}{2} \vspace{.05in}\\
\phantom{\frac{1}{2}} & \frac{1}{4}, & 1-i\frac{p}{2}, & 1+i\frac{p}{2}, &
1-\frac{p}{2}, & 1+\frac{p}{2} \end{array}
\Big| \; -1 \right]\\
\equiv \sum^{\frac{p-1}{2}}_{k=0}
(4k+1) {\bin{-\frac{1}{2}}{k}}^5  \pmod {p^4} \; .
\end{multline}

\noindent
Next we examine the right hand side of (\ref{Whipple}). By (\ref{GHS_def}),
\begin{multline}\label{Whipple_RHS}
\frac{\biggf{1-\frac{p}{2}}\biggf{1+\frac{p}{2}}}
{\biggf{\frac{3}{2}}\biggf{\frac{1}{2}}}
\; {_3F_2} \left[ \begin{array}{ccc} \frac{1}{2} & \frac{1}{2} + \frac{p}{2}, & \frac{1}{2} - \frac{p}{2} \vspace{.05in}\\
\phantom{\frac{1}{2}} & 1-i\frac{p}{2}, & 1+i\frac{p}{2} \end{array}
\Big| \; 1 \right]\\
=\frac{\biggf{1-\frac{p}{2}}\biggf{1+\frac{p}{2}}}
{\biggf{\frac{3}{2}}\biggf{\frac{1}{2}}}
\sum^{\frac{p-1}{2}}_{k=0}
\frac{\Bigph{\frac{1}{2}}{k}
\Bigph{\frac{1}{2}+\frac{p}{2}}{k} \Bigph{\frac{1}{2}-\frac{p}{2}}{k}}
{\Bigph{1-i\frac{p}{2}}{k} \Bigph{1+i\frac{p}{2}}{k}}
\; \frac{1}{{k}!} \; .
\end{multline}

\noindent
Now, via (\ref{reflect}) and the fact that $\gf{x+1}=x\gf{x}$ and $\gf{\frac{1}{2}}=\sqrt{\pi}$, we have
\begin{align}\label{RHS_Sub1}
\frac{\biggf{1-\frac{p}{2}}\biggf{1+\frac{p}{2}}}
{\biggf{\frac{3}{2}}\biggf{\frac{1}{2}}}
&=\frac{\biggf{1-\frac{p}{2}} \bigl({\frac{p}{2}}\bigr) \biggf{\frac{p}{2}}}
{\biggf{\frac{3}{2}}\biggf{\frac{1}{2}}}\\
%\notag &=\frac{\Bigl({\frac{p}{2}}\Bigr) \Bigl({ \frac{ \pi}{\sin{ \left({ \frac{p}{2} \pi}\right) }} }\Bigr)}
%{\frac{\sqrt{\pi}}{2} \sqrt{\pi}}\\
\notag &=\frac{p}{\sin{ \left({ \frac{p}{2} \pi}\right) }}\\
\notag &=\phi_p(-1)\:p \; .
\end{align}
Also, we have
\begin{align}\label{RHS_Sub2}
\frac
{\Bigph{\frac{1}{2}+\frac{p}{2}}{k} \Bigph{\frac{1}{2}-\frac{p}{2}}{k}} 
{\Bigph{1-i\frac{p}{2}}{k}  \Bigph{1+i\frac{p}{2}}{k}}
\equiv {\bin{-\frac{1}{2}}{k}}^2  \pmod {p^2} \;.
\end{align}

\noindent
Using (\ref{LHS_Sub1}) and substituting (\ref{RHS_Sub1}), (\ref{RHS_Sub2}) into (\ref{Whipple_RHS}), we get
\begin{multline}\label{RHS_Final}
\frac{\Biggf{1-\frac{p}{2}}\Biggf{1+\frac{p}{2}}}
{\Biggf{\frac{3}{2}}\Biggf{\frac{1}{2}}}
\; {_3F_2} \left[ \begin{array}{ccc} \frac{1}{2} & \frac{1}{2} + \frac{p}{2}, & \frac{1}{2} - \frac{p}{2} \vspace{.05in}\\
\phantom{\frac{1}{2}} & 1-i\frac{p}{2}, & 1+i\frac{p}{2} \end{array}
\Big| \; 1 \right]\\
\equiv \phi_p(-1) \: p
\left[{ \sum^{\frac{p-1}{2}}_{j=0} {(-1)}^j {\bin{-\frac{1}{2}}{j}}^3 }\right] \pmod {p^3} \; .
\end{multline}

\noindent
Finally, combining (\ref{Whipple}), (\ref{LHS_Final}) and (\ref{RHS_Final}) yields (\ref{lemma3_2}) 
and hence the result follows. 
\end{proof}
%%%%%%%%%%%%%%% End Lemma 3 - Z(p,1,2) %%%%%%%%%%%%%%%%%%%%%%

\section*{acknowledgements}
The first author would like to thank the UCD Ad Astra Research Scholarship programme for its financial support. The second author thanks the Institut des Hautes {\'E}tudes
Scientifiques for their hospitality and support during the preparation
of this paper.

\end{document}